\documentclass[12pt, twoside]{article}
\usepackage{amsmath,amsthm,dsfont}
\usepackage{amsfonts}
\usepackage{amssymb,latexsym}
\usepackage{enumerate}
\newtheorem{theorem}{Theorem}[section]

\newtheorem{lemma}[theorem]{Lemma}

\newtheorem{exam}[theorem]{Example}
\newtheorem{corollary}[theorem]{Corollary}
\def\n{\noindent}
\def\fr{\frac}
\def\Om{\Omega}
\def\loc{\text{loc}}

\def\de{\delta}
\def\ve{\varepsilon}
\def\va{\varphi}
\def\la{\lambda}

\begin{document}
\setlength{\baselineskip}{18truept}
\pagestyle{myheadings}
\markboth{ }{Approximation in holomorphic Bergman spaces}
\title {Approximation in holomorphic Bergman spaces}
\author{\;\;\; \;\;\; 
Nguyen Van Phu\\
 Faculty of Natural Sciences, Electric Power University\\ Hanoi, Vietnam\\E-mail: phunv@epu.edu.vn}
\date{}
\maketitle 

\renewcommand{\thefootnote}{}



\renewcommand{\thefootnote}{\arabic{footnote}}
\setcounter{footnote}{0}

\begin{abstract}
Let $\varphi, \psi $ be negative plurisubharmonic functions in a pseudoconvex domain  $\Omega$ in $\mathbb{C}^{n}$ and $f$ be a holomorphic function  belonging to $L^{2}(\Omega, \varphi)$. We give a sufficient condition so that $f$ can be approximated in $L^2$ norm by elements in $L^{2}(\Omega^{'}, \psi)$ where $\Omega^{'}$ is a relative compact subset of $\Omega.$
\end{abstract}
\section{Introduction}
Let $\Omega$ be a domain in $\mathbb{C}^{n}.$ By $PSH^-(\Omega),$ we denote the set of negative plurisubharmonic functions on $\Omega$ and by $Hol(\Omega),$ we denote the set of holomorphic functions on $\Omega.$ If $\phi: \Omega\mapsto [-\infty,+\infty)$ is an upper semi-continuous function on $\Omega,$ we define
$$L^2(\Omega,\phi)=\{f:\Omega\to\mathbb{C}: ||f||^{2}_{L^2(\Omega,\phi)}=\int_{\Omega}|f|^2e^{-\phi}dV_{2n}<\infty\}.$$
In the case where $\phi\equiv 0,$ for convinience, we use the notation $L^2(\Omega)$ instead of $L^2(\Omega,0).$\\
\n  Let $\varphi,\psi\in PSH^-(\Omega)$ and $f$ be a holomorphic function which belongs to $L^2(\Omega,\varphi)$.
A question arises:  when $\psi$ is close to $\varphi$ in the $L^1_{\loc}$ norm, can  we can choose a holomorphic function $g$ belonging to the class $L^2(\Omega,\psi)$ that is close to $f$ in the $L^2_{\loc}$ norm.\\
To this purpose, we first provide a sufficient condition to make this happen.
\begin{lemma}\label{key}
Let $\Omega \subset \mathbb C^n$ be a bounded pseudoconvex domain and 
$\psi \in$ PSH$^{-}(\Omega).$ Then for every $f\in Hol(\Omega) \cap L^2 (\Om)$ 
and $t>1$ we can find $g \in Hol (\Om) \cap L^2 (\Omega, \psi)$
such that
\begin{equation} \label{eq0}
\fr1{2}\Vert g-f\Vert_{L^2 (\Om)}^2 \le \int\limits_{\{\psi<-t\}} |f|^2 dV_{2n}+
144t^2 e^{t+1}\int\limits_{\{-(t+1)<\psi<-t\}} |f|^2dV_{2n}.
\end{equation}
\end{lemma}
\n The main result of our note give a sufficient condition so that a holomorphic function $f\in L^2(\Omega,\varphi)$ can be approximated in $L^2$ norm by elements in $L^{2}(\Omega^{'}, \psi)$ where $\Omega^{'}$ is a relative compact subset of $\Omega.$
\begin{theorem} \label{main}
Let $\Omega \subset \mathbb C^n$ be a bounded pseudoconvex domain and 
$\va \in$ PSH$^{-}(\Omega).$  Let $K$ be a compact subset of $\Om$ and
$f\in Hol(\Omega) \cap L^2 (\Om, \va)$.
Then for every $\ve>0$ there exists $\de>0$ such that for all $\psi \in PSH^{-} (\Om)$	
satisfying $\Vert \va-\psi\Vert_{L_{\loc}^1 (\Om)}<\de$
we can find  $g \in Hol (\Om) \cap L^2 (\Omega^{'}, \psi)$
such that $\Vert g-f\Vert_{L^2 (K)}^2<\ve,$ where $K\subset\Omega^{'}\Subset\Omega.$
\end{theorem}
\n Theorem \ref{main} implies a result concerning stability of the multiplier ideal sheaf.
\begin{corollary} \label{coro}
	Let $\Omega \subset \mathbb C^n$ be a bounded pseudoconvex domain. Assume that
	 $\{\varphi_j\}_{j\geq1}\subset PSH(\Omega),$  $\varphi\in PSH(\Omega)$ be such that $\varphi_j\leq \varphi$ for all $j\geq 1$ and $\varphi_j\to\varphi$ in $L^1_{\loc}(\Omega).$ Then for $\Omega^{'}\Subset\Omega$ there exist $j_{0}\geq 1$ such that $\mathcal{I}(\varphi_j)=\mathcal{I}(\varphi) $ on $\Omega^{'}$ for all $j\geq j_0.$
\end{corollary}
Our note is organized as follows. In Section 2, we recall some auxiliary facts about $L^{2}$- estimate for the $\overline{\partial}-$operator and Demailly's strong openness conjecture. In Section 3, we give the proof of the main Theorem of this paper. Finally, in Section 4, we give some examples to illustrate the main theorem.
{\bf Acknowledgements.}  
The author is supported by the grant B2025-CTT-10 from Ministry of Education and Training, Vietnam.
\section{Preliminaries}
\subsection{$L^{2}$- estimate for the $\overline{\partial}-$operator}
In this Subsection, we recall some results about $L^{2}$- estimate for the $\overline{\partial}-$operator. We refer the readers to  $\cite{BL16}$, $\cite{Bl05}$, $\cite{De97}$, $\cite{Ho91}$ for further information.\\
Let $\Omega$  be a pseudoconvex domain in $\mathbb{C}^{n}$. We denote by $L^2_{\loc}(\Omega)$ the set of functions in $\Omega$ which are locally square integrable with respect to the Lebesgue measure and by $L^{2}_{\loc,(0,1)}(\Omega)$ we denote the set of forms of type $(0,1)$ with coefficients in $L^2_{\loc}(\Omega).$ Assume that 
$$\alpha = \sum_{j=1}^{n}\alpha_{j}d\overline{z}_{j}\in L^{2}_{\loc,(0,1)}(\Omega)$$ 
satisfies $\overline{\partial}\alpha = 0.$ We are interested in solving 
\begin{equation}\label{ep1}
\overline{\partial}u=\alpha
\end{equation}
with an $L^{2}-$estimate. Such solutions are very useful in constructing new holomorphic functions because $\overline{\partial}v=0$ implies that  $v$ is holomorphic.
The classical one is due to L. H\"ormander $\cite{Ho91}$: For a smooth, strongly plurisubharmonic function $\varphi$ in $\Omega,$ one can find a solution of ($\ref{ep1}$) satisfying 
\begin{equation}\label{ep2}
\int_{\Omega}|u|^{2}e^{-\varphi}d\lambda\leq \int_{\Omega}|\alpha|^{2}_{i\partial\overline{\partial}\varphi}e^{-\varphi}d\lambda
\end{equation}
where
$$|\alpha|^{2}_{i\partial\overline{\partial}\varphi}=\sum_{j,k}\varphi^{j,\overline{k}}\overline{\alpha}_{j}\alpha_{k}$$
 is the length of $\alpha$ with respect to the K\"ahler metric $i\partial\overline{\partial}\varphi$ with potential $\varphi.$ Here $(\varphi^{j,\overline{k}})$ is the inverse transposed of the complex Hessian $(\dfrac{\partial^{2}\varphi}{\partial z_{j}\partial\overline{z}_{k}}).$
It was observed in Theorem 3.1 in $\cite{Bl05}$ that the H\"ormander estimate ($\ref{ep2}$) also holds for arbitrary plurisubharmonic function $\varphi$  if instead of $ |\alpha|^{2}_{i\partial\overline{\partial}\varphi}$ one should take  any nonnegative $H\in L^{\infty}_{loc}(\Omega)$ with
$$i\overline{\alpha}\wedge\alpha\leq H i\partial\overline{\partial}\psi,$$ where $\psi$ is a plurisubharmonic function such that $-e^{-\psi}$ is also plurisubharmonic.
Then we can find $u\in L_{\loc}^2(\Omega)$ with $\overline{\partial}u=\alpha$ and such that 
$$\int_{\Omega}|u|^{2}e^{-\varphi}d\lambda\leq 16\int_{\Omega}He^{-\varphi}d\lambda.$$
\subsection{ Demailly's strong openness conjecture}
The singularities of the plurisubharmonic function plays an important role in several complex variables and complex geometry. Various properties about the singularities of the plurisubharmonic function have been discussed ( e.g., see  \cite{Nadel}, \cite{De01}, \cite{DP03}, \cite{H14}, \cite{GZ}).
In this Subsection, we recall the strong openness conjecture about multiplier ideal sheaves for plurisubharmonic functions, which was posed by Demailly in \cite{DP03}. \\
\n Denote by $\mathcal{O}_{\mathbb{C}^n,z}$ the ring of germs of holomorphic functions at $z$. Let $\varphi\in PSH^-(\Omega).$ Following \cite{Nadel}, one can define the multiplier ideal sheaf $\mathcal{I}({\varphi})$  to be the sheaf of germs of holomorphic function $f\in \mathcal{O}_{\mathbb{C}^n,z}$ such that $\int_{U}|f|^{2}e^{-\varphi}dV_{2n}<+\infty$ on some neighborhood $U$ of $z$. \\
{\it Strong openness conjecture}. Let $\varphi$ be a plurisubharmonic function on $\Omega.$ Then $$\mathcal{I}(\varphi)=\mathcal{I}_{+}(\varphi):=\cup_{\varepsilon> 0}\mathcal{I}\big((1+\varepsilon)\varphi\big).$$
Recently,  Q. Guan and X. Zhou \cite{GZ} and P. H.  Hiep  \cite{H14}  proved strong openness conjecture is true.
 In fact, this statement 
is equivalent to the following theorem (see Theorem 1.1 in \cite{GZ}).
\begin{theorem} \label{bd1}
Let $\varphi$ be a negative plurisubharmonic function on the unit polydisk $\Delta^n(0,1)\subset\mathbb{C}^n.$ Suppose that $f\in L^2(\Delta^n(0,1),\varphi)$ is a holomorphic function on $\Delta^n(0,1).$ Then for some $r\in(0,1),$ there exists a constant $\lambda>1$ such that $f\in L^2\big(\Delta^n(0,r), \lambda \varphi\big).$
\end{theorem}
\section{Proof of Main results}
\begin{proof} (of Lemma \ref{key})
The idea to construct such a holomorphic function $g$ is taken from the comment before Theorem 3.3 in $\cite{Bl05}.$ Firstly, we put
$$\phi=-\log(-\psi).$$
 It follows that  $$-e^{-\phi}=\psi\in PSH^{-}(\Omega).$$
 This implies
 $$\partial\overline{\partial}(-e^{-\phi})\geq 0.$$
 Moreover, we have
 \begin{align*}
 \partial\overline{\partial}(-e^{-\phi})&=\partial(e^{-\phi}\overline{\partial}\phi)\\
 &=-e^{-\phi}\partial\phi\wedge \overline{\partial}\phi+ e^{-\phi} \partial\overline{\partial}\phi\\
 &=e^{-\phi}[\partial\overline{\partial}\phi-\partial\phi\wedge \overline{\partial}\phi].
 \end{align*}
So we obtain
$$\partial\overline{\partial}\phi\geq\partial\phi\wedge \overline{\partial}\phi .$$
Now for $a>0$ we define the following function 
\begin{equation}\label{e1}
	\chi(t)=
\begin{cases}
	1&\,\, \text{if}\,\, t<a\\
	0 &\,\, \text{if}\,\, t>\log(e^a+1)\\
	\alpha (t-\log(e^a+1))^3+\beta (t-\log(e^a+1))^2&\,\, \text{if}\,\, a<t<\log(e^a+1),
\end{cases}
\end{equation}
where 
$$\alpha:=-\frac2{(a-\log(e^a+1))^3}, \beta:= \frac3{(a-\log(e^a+1))^2}.$$
It is easy to check that, $\chi$ is $\mathcal C^1$ smooth and $0 \le \chi \le 1$ on $\mathbb R$
and 
\begin{equation} \label{eq6}
-3e^a \le -\frac3{2(\log (1+\frac1{e^a}))}=-\frac3{2(\log (e^a+1)-a)} \le \chi'(t) \le 0 \ \forall t \in \mathbb R.
\end{equation}
Here the first inequality follows from the fact that 
$$\log \Big (1+\fr1{x} \Big ) \ge \fr{1}{2x}, \ \forall x=e^a>1.$$
Notice that
$$\chi(\log(-\psi))= 
\begin{cases}
	1&\,\, \text{if}\,\, \psi>-e^a\\
	0 &\,\, \text{if}\,\, \psi<-(e^a+1)\\
	\in [0;1 ]&\,\, \text{if}\,\, -(e^a+1)<\psi<-e^a.
\end{cases}$$
We shall use Theorem 3.1 in $\cite{Bl05}$ for $\alpha$ given by
$$\alpha= -f\overline{\partial}\big(\chi(\log(-\psi))\big).$$ 
Indeed, we have 
\begin{align*}
 \alpha&= -\overline{\partial}\big(\chi(-\phi)f\big)\\
&=f\chi^{'} \big(\log(-\psi)\big).\overline{\partial}\phi\\
&=f\chi^{'} \big(\log(-\psi)\big).\dfrac{\overline{\partial}\psi}{|\psi|}.
\end{align*}
This implies
\begin{align*}
	i\overline{\alpha}\wedge \alpha=&i|f|^{2}| \chi^{'} \big(\log(-\psi)\big)|^2\dfrac{\partial\psi\wedge \overline{\partial}\psi}{|\psi|^{2}}\\
	&=i|f|^{2}| \chi^{'} \big(\log(-\psi)\big)|^{2}\partial\phi\wedge \overline{\partial}\phi\\
	&\leq i|f|^{2}| \chi^{'} \big(\log(-\psi)\big)|^{2}\partial\overline{\partial}\phi\\
	&\leq iH\partial\overline{\partial}\phi	
\end{align*}
with $$ H=|f\chi^{'} \big(\log(-\psi)\big)|^{2}.$$
Let $u$  be solution to  $$\overline{\partial}u=\alpha$$ which is minimal in the $L^{2}(\Omega, \psi)$ norm. 
Now, we put $$g=u + \chi\big(\log(-\psi)\big)f$$
then $\bar \partial g=0,$ so $g\in Hol(\Omega)$. Observe that
$$\begin{aligned} 
\int\limits_{\Omega}| \chi\big(\log(-\psi)\big)|^{2}|f|^{2}e^{-\psi}dV_{2n}
&\leq \int\limits_{\{\psi \ge -(e^a+1)\}}|f|^{2}e^{-\psi}dV_{2n}\\
&\leq e^{e^a+1} \int\limits_{\Om}|f|^{2}dV_{2n}<\infty.
\end{aligned}$$
Since $u \in L^{2}(\Omega, \psi)$, we conclude that $g \in L^{2}(\Omega, \psi).$

Next we write
\begin{align}\label{eq3}
	\begin{split}
	\fr1{2} \Vert g-f\Vert_{L^2 (\Om)}^2&=\fr1{2}\int\limits_{\Omega}|g-f|^{2}dV_{2n}\\
	&=\fr1{2}\int\limits_{\Omega}|u + \chi\big(\log(-\psi)\big)f-f|^{2}dV_{2n}\\
	&\leq  \int\limits_{\Omega}|u|^2dV_{2n} + \int\limits_{\Omega} (1-\chi\big(\log(-\psi)\big)^2|f|^{2}dV_{2n}.	
	\end{split}
\end{align}
Next
\begin{align*}
\int_{\Omega} (1-\chi\big(\log(-\psi)\big))^2|f|^{2}dV_{2n}
&= \int\limits_{\{\psi\leq -e^a\}} (1-\chi (\log(-\psi)))^2|f|^{2}dV_{2n}\\
&\leq \int\limits_{\{\psi\leq -e^a\}} |f|^{2}dV_{2n}.
\end{align*}
Finally, using Theorem 3.1 in \cite{Bl05} we obtain
$$\begin{aligned}
 \int\limits_{\Omega}|u|^2dV_{2n}
&\leq \int\limits_{\Omega}|u|^2e^{-\psi}dV_{2n}\\
 &\leq 16\int\limits_{\Omega}He^{-\psi}dV_{2n}\\
 &=16\int_{\Omega}|f\chi^{'} \big(\log(-\psi)\big)|^{2}e^{-\psi}dV_{2n}\\
  &\leq 144e^{2a}\int\limits_{\{-(e^a+1)<\psi<-e^a\}} |f|^2e^{-\psi}dV_{2n}.
  \end{aligned}$$
Putting all this together and regarding $t=e^a,$ we complete the proof of Lemma \ref{key}.
\end{proof}
Now we turn to the proof of the main theorem.
\begin{proof}
Let $\rho$ be a smooth strictly plurisubharmonic exhaution function for $\Om$ such that $\rho<-1$ on $K$.
Define $$\Om{''}:=\{\rho<0\}, \Om':=\{\rho<-1\}.$$ Then $\Om' \Subset \Om{''}$ are relatively compact pseudoconvex open subsets of $\Om$ saitisfying $K \subset \Om'.$ 	
We split the proof into two steps.

\n 
{\it Step 1.} Given $\ve>0$ we will show that there exists $\de>0$ such that for all $\psi \in PSH^{-} (\Om{''})$	
satisfying $\Vert \va-\psi\Vert_{L^1 (\Om{''})}<\de$
we can find  $g \in Hol (\Om{''}) \cap L^2 (\Omega^{''}, \psi)$
with $\Vert g-f\Vert_{L^2 (\Om{''})}^2<\ve.$ 	
	
To this purpose, we may assume $\vert f \vert \le 1$ on $\Om^{''}.$
Since $f \in Hol(\Om) \cap L^2(\Om, \va),$ by the strong openess property (Theorem \ref{bd1}) we can find a constant $\la>1$ such that
 $$\int\limits_{\Om'} |f|^2 e^{-\la \va} dV_{2n}<\infty.$$
Since the pluripolar set $\{\va=-\infty\}$ has the Lebesgue measure $0,$
we may find $t_0>1$ so large such that
$$\int\limits_{\{\va<-t_0+1\} \cap \Om^{''}} |f|^2e^{-\la \va}dV_{2n}<\ve.$$
This implies that
$$  \int\limits_{\{\va<-t_0+1\} \cap \Om^{''}} |f|^2 dV_{2n} \le \fr{\ve}{e^{\la (t_0-1)}}.$$
Let $\psi \in PSH^{-} (\Om^{''})$ be such that 
$$\Vert \va-\psi\Vert_{L^1 (\Om^{''})}<\de:= \fr{\ve}{4(144t_0^2 e^{t_0+1}+1)}.$$
Using Lemma \ref{key} 
we may find $g \in Hol(\Om^{''}) \cap L^2 (\Om^{''}, \psi)$ satisfying (\ref{eq0}).
Then we have
$$\begin{aligned}
\int\limits_{\{-(t_0+1)<\psi<-t_0\} \cap \Om^{''}} |f|^2dV_{2n} & \le \int\limits_{\{\psi<-t_0\} \cap \Om^{''}} |f|^2dV_{2n}\\ &
\le \int\limits_{\{\va<-t_0+1\} \cap \Om^{''}} |f|^2dV_{2n}+\int\limits_{\{\va \ge \psi+1\} \cap \Om^{''}} |f|^2 dV_{2n}\\
&\le \fr{\ve}{e^{\la (t_0-1)}} +\int\limits_{\Om^{''}} |\va-\psi|dV_{2n}\\
&<\de+\de=2\de,
\end{aligned}$$	
provided that $t_0$ is taken so large  that
$e^{\la (t_0-1)}>4(144t_0^2 e^{t_0+1}+1).$

Combining all these esttimates we get 
$$\fr1{2}\Vert g-f\Vert^2_{L^2 (\Om^{''})} \le 2\de [1+144t_0^2e^{t_0+1}].$$
After rearranging the above inequality and plugging the value of $\de$ we obtain
the desired estimate $\Vert g-f\Vert_{L^2 (\Om^{''})}^2<\ve.$ 

\n
{\sl Step 2.} Completion of the proof. By Step 1 we may find $\de'>0$ such that if $\psi \in PSH^{-} (\Om)$ and
satisifes $\Vert \psi-\va\Vert_{L^1_{\loc} (\Om)}<\de'$
then there exists $g \in Hol(\Om^{''}) \cap L^2 (\Om^{''}, \psi)$ with 
\begin{equation} \label{eq11}
\Vert g-f\Vert_{L^2 (\Om^{''})}^2<\fr{\ve}4.
\end{equation}
Let $0 \le \chi \le 1$ be a smooth cut-off function with support in $\Om^{''}$ such that $\chi=1$ on a neighbourhood of $\overline{\Om'}$. Then an easy calculation yields that
$|\bar \partial (\chi g)| \in L^2 (\Om, \psi).$
Next we solve the equation 
$\bar \partial u_j=\bar \partial (\chi g)$ with $L^2$-H\"ormander estimate for the weight 
$e^{-j\max \{\rho+1,0\}-\psi}$ to get a solution $u_j$ satisfying the following estimate
$$\int_{\Om} |u_j|^2 e^{-j\max \{\rho+1,0\}-\psi}dV_{2n} \le C \int_{\Om} |\bar \partial (\chi g)|^2 e^{-j\max \{\rho+1,0\}-\psi}dV_{2n}<\infty,$$
where $C>0$ depends only on $n$ and the diameter of $\Om.$ 
Thus $g_j:=\chi g-u_j$ is holomorphic on $\Om.$ Since the hypersurface $\{\rho=-1\}$ has Lebesgue measure $0$, 
by the choice of $\chi$ we infer
that
$$\int\limits_{\Om} |\bar \partial (\chi g)|^2 e^{-j\max \{\rho+1,0\}-\psi}dV_{2n}
\le \int\limits_{\{\rho>-1\} \cap \Om^{''}}  |\bar \partial (\chi g)|^2 e^{-j\max \{\rho+1,0\}-\psi}dV_{2n}.$$
Since the right terms converge to $0$ as $j \to \infty$ by Lebesgue's dominated convergence theorem we can find $j$ so large such that
$$\begin{aligned}
\int\limits_{\Om^{'}} |u_j|^2 dV_{2n} &\le \int\limits_{\Om^{'}} |u_j|^2 e^{-\psi} dV_{2n}\\
&\le \int\limits_{\Om} |\bar \partial (\chi g)|^2 e^{-j\max \{\rho+1,0\}-\psi}dV_{2n}\\
&\le \fr{\ve}4.
\end{aligned}$$  
It implies that $g_j:=\chi g-u_j \in L^2(\Om^{'},\psi)$ and note that since $\chi=1$ on $K$, we obtain
$\Vert g_j-g\Vert^2_{L^2 (K)} \le  \ve/4$. Combining this estimate and (\ref{eq11}) we obtain
	$$\Vert g_j-f\Vert^2_{L^2 (K)} \le 2 [\Vert g_j-g\Vert^2_{L^2 (K)}+\Vert g-f\Vert^2_{L^2 (K)}]<\ve.$$
The proof is then complete.
\end{proof}
\n Next we give a proof of Corollary \ref{coro}
\begin{proof}
It follows from $\varphi_j\leq \varphi$ that $\mathcal{I}(\varphi_j)\subset \mathcal{I}(\varphi).$ The reverse inclusion follows directly from the main theorem.
\end{proof}
\section{Examples}
In this Section, taking idea from Remark 1.3 in \cite{H14}, we give examples to illustrate the main theorem.
\begin{exam}\label{exam1}
\n  Assume that $\Omega = \Delta^2(0,r)=\Delta(0,r) \times\Delta(0,r)\subset \mathbb{C}^{2}$	
and $0<\varepsilon<1,$ let us choose $\varphi_{j}(z)= (4-\varepsilon)\ln|z_1+\frac{z_{2}}{j}|, \varphi(z)=(4-\varepsilon)\ln|z_1|$ and $f(z) = z_1.$ Then we have $\varphi_{j}(z)\to\varphi(z)$ in $L^1_{\loc}(\Omega)$ as $j\to+\infty$. \\
Using Fubini's theorem, we see that
\begin{align*}\int_{\Delta^2(0,r)} e^{-\varphi}|f|^2dV_{4} &=\int_{\Delta^2(0,r)} e^{-2(2-\frac{\varepsilon}{2})\ln|z_1|}|z_1|^2dV_{4} = \int_{\Delta^2(0,r)} \dfrac{1}{|z_1|^{2-\varepsilon}}dV_{4}\\
	&= V_2(\Delta(0,r))\int_{\Delta(0,r)} \dfrac{1}{|z_1|^{2-\varepsilon}}dV_{2}.
	\end{align*} 
Using change of variables $z_1=\rho e^{i\theta}$ with $\theta\in[0,2\pi]$ and $\rho\in[0,r],$ we obtain
\begin{equation}\label{e4.1}
	\int_{\Delta(0,r)} \dfrac{1}{|z_1|^{2-\varepsilon}}dV_{2}=\int_{0}^{2\pi}d\theta\int_{0}^{r}\frac{\rho d\rho}{\rho^{2-\varepsilon}}=2\pi\int_{0}^{r}\frac{d\rho}{\rho^{1-\varepsilon}}<+\infty.
	\end{equation}
That means we have $f\in L^2(\Omega,\varphi).$ Obviously, $f\in Hol(\Omega).$\\
Firstly, we show that in general, we cannot choose $g=f$ in the Main Theorem. Indeed, for $0<\rho<r$, we have $ \Delta^2(0,\rho)\Subset\Omega.$ Therefore,
\begin{align*}
	\int_{\Delta^2(0,\rho)} e^{-\varphi_{j}}|f|^2dV_{4}
	&=\int_{\Delta^2(0,\rho)}e^{-2(2-\frac{\varepsilon}{2})\ln|z_1+\frac{z_{2}}{j}|}|z_1|^2dV_{4} \\
	&=\int_{\Delta^2(0,\rho)}\dfrac{1}{|z_1+\dfrac{z_{2}}{j}|^{4-\varepsilon}.|z_{1}|^{-2}}dV_{4}.
\end{align*}
Using change of variables $$\begin{cases}
	w_1&=z_1+\frac{z_{2}}{j}\\
	w_2&=z_2.
	\end{cases}$$
Then Jacobi determinant $J=1.$ Thus, using Fubini's theorem and repeating the argument as in (\ref{e4.1}), we have
\begin{align*}
	\int_{\Delta^2(0,\rho)} e^{-\varphi_{j}}|f|^2dV_{4}
	&=\int_{\Delta(0,\rho_1)\times \Delta(0,\rho) } \dfrac{1}{|w_1|^{4-\varepsilon}.|w_{2}|^{-2}}dV_{4}\\
	&=\int_{\Delta(0,\rho_1}\dfrac{1}{|w_1|^{4-\varepsilon}}dV_{2}. \int_{ \Delta(0,\rho)} \dfrac{1}{|w_{2}|^{-2}}dV_{2} \\
	&=+\infty.
\end{align*}	
This means that, $f\notin L^2(\Delta^2(0,\rho),\varphi_j)$ and we cannot choose $g=f$.\\
Secondly,  we give an example of a holomorphic function $g$ satisfies  Theorem \ref{main}.
Indeed, we can choose $g(z)=z_{1}+\dfrac{z_{2}}{j}.$ Obviously, we have $g\in Hol(\Omega)$ and $||g-f||_{L^2(K)}\to 0$ with $j$ large enough. Moreover, repeating the above argument, we also have
\begin{align*}
	\int_{\Delta^2(0,r)} e^{-\varphi_{j}}|g|^2dV_{4}
	&=\int_{\Delta^2(0,r)}e^{-2(2-\frac{\varepsilon}{2})\ln|z_1+\frac{z_{2}}{j}|}|z_1+\dfrac{z_{2}}{j}|^2dV_{4} \\
	&=\int_{\Delta^2(0,r)}\dfrac{1}{|z_1+\dfrac{z_{2}}{j}|^{2-\varepsilon}}dV_{4}\\
	&<+\infty.
\end{align*}
It means that $g\in L^2(\Delta^2(0,r),\psi)$ and $g$ satisfies the Theorem $\ref{main}.$
\end{exam}
\begin{exam} In this Example, we show that one cannot remove condition $\varphi_j\leq \varphi$ (which implies that $ \mathcal{I}({\varphi_{j}})\subset \mathcal{I}({\varphi})$)  in the Corollary $\ref{coro}.$\\
	\n  Indeed, for simplicity,  we concern multiplier ideal sheaves that are subset of $\mathcal{O}_{\mathbb{C}^n,0}.$ Let $\Omega,\varphi_{j},\varphi$ as in Example \ref{exam1}. Then we have $f(z)=z_{1}\in\mathcal{I}(\varphi)$ and $g(z)=z_{1}+\dfrac{z_{2}}{j}\in \mathcal{I}(\varphi_{j}).$ Assume that if
	$\mathcal{I}({\varphi_{j}})= \mathcal{I}({\varphi}),$ then we have 
	$g(z)\in \mathcal{I}(\varphi).$ Thus, we obtain $z_{2}\in \mathcal{I}(\varphi).$ This is not true since for all $0<\rho<r$, we have
	\begin{align*}
		\int_{\Delta^2(0,\rho)} e^{-\varphi}|z_{2}|^2dV_{4}
		&=\int_{\Delta^2(0,\rho)}e^{-2(2-\frac{\varepsilon}{2})\ln|z_1|}|z_2|^2dV_{4} \\
		&=\int_{\Delta^2(0,\rho)}\dfrac{1}{|z_1|^{4-\varepsilon}|z_{2}|^{-2}}dV_{4}\\
		&=+\infty.
	\end{align*}
\end{exam}
\section*{Declarations}
\subsection*{Data Availability}
This declaration is not applicable.
\subsection*{Competing interests}
The authors have no conflicts of interest to declare that are relevant to the content of this article.
\subsection*{Funding }
The paper is supported by the grant B2025-CTT-10 from Ministry of Education and Training, Vietnam.
\subsection*{Ethics and Consent to Participate declarations}
This declaration is not applicable.

\end{document}